\documentclass{amsart}

\usepackage{eufrak}

\usepackage{amssymb,amsmath, amscd, mathrsfs, yfonts, bbm, graphics}

\title[multiplicative
boolean cumulants] {A new proof for the multiplicative property of
the boolean cumulants with applications to operator-valued case}

\author{Mihai Popa}
\address{${}^1$Indiana University at Bloomington,
 Department of Mathematics, Rawles Hall,
 931 E 3rd St, Bloomington, IN 47405}
\email{mipopa@indiana.edu}
\address{${}^2$Institute of Mathematics, Romanian Academy,
 P.O.Box 1-764, Bucharest, RO-70700, Romania}
\DeclareMathAlphabet{\mathpzc}{OT1}{pzc}{m}{it}
\newtheorem{claim}{}[section]
\newtheorem{defn}[claim]{Definition}

\newtheorem{remark}[claim]{Remark}
\newtheorem{prop}[claim]{Proposition}
\newtheorem{cor}[claim]{Corollary}

\newcommand{\cA}{\mathcal{A}}
\newcommand{\cC}{\mathcal{C}}

\newcommand{\wb}{\widetilde{b}}
\newcommand{\wm}{\widetilde{m}}
\newcommand{\I}{\mathcal{I}}

\newcommand{\ou}{\overrightarrow{u}}
\newcommand{\ov}{\overrightarrow{v}}
\newcommand{\ow}{\overrightarrow{w}}

\newcommand{\wB}{\widetilde{B}}
\newcommand{\wM}{\widetilde{M}}
\newcommand{\gB}{\mathfrak{B}}

\newcommand{\gA}{\mathfrak{A}}

\newcommand{\lra}{\longrightarrow}

\newcommand{\p}{^\prime}

\newcommand{\mul}{Mul[[\gB]]}

\usepackage{amssymb}
\input xy
\xyoption{all}

\begin{document}
 \maketitle
\bibliographystyle{alpha}
\begin{abstract}
 The paper presents several combinatorial properties of the boolean
 cumulants. A corollary is a new proof of the multiplicative
 property of the boolean cumulant series that can be easily
 adapted for the case of boolean independence with amalgamation over
 an algebra.


\end{abstract}

\section{Introduction}
 Boolean probability theory, according to \cite{Schurmann} is
 one of the three symmetric universal
probability theories (the other two being classical and free
probability theories). It have been in the literature at least since
early '70's (\cite{vWa}) with various developments, from stochastic
differential equations (\cite{Ter}) to measure theory
(\cite{speiwa}, \cite{ber}) and Appell polynomials (\cite{anshel}).

Two important notions in free probability theory are the $R$ and $S$
transforms. For $X$ a non-commutative random variable of finite
moments, $R_X(z)$ and $S_X(z)$ are power series that encode the
information from the moment generating function of $X$ and have the
following additive, respectively multiplicative, property: if $X,Y$
are free, then $R_{X+Y}=R_X=R_Y$ and $S_{XY}=S_XS_Y$.

The notion of boolean cumulants has been in the literature, in
various contexts, at least since early '70's (\cite{vWa},
\cite{Ter}, \cite{Hud}). Their generating function, that we will
call the $B$-transform, has an additive property similar to the
$R$-transform (\cite{speiwa}): if $X$ and $Y$ are boolean
independent, then $B_{X+Y}=B_X+B_Y$.

In \cite{Uwe1}, U. Franz remarked an interesting multiplicative
property of the $B$-transform, which gives a valuable tool for
studying the multiplicative boolean convolution. For $X,Y$ boolean
independent, one has that
\begin{equation}\label{multipli}
B_{(1+X)(1+Y)}=B_{(1+X)}B_{(1+Y)}. \end{equation}

In literature there are two proves of \ref{multipli}. The original
proof, in \cite{Uwe1}, uses algebraic properties of the resolvent
function. The proof in \cite{ber}, in the spirit of D. Voiculescu's
proof of the property for the $R$ and $S$-transforms, uses a
Fock-type model and properties of the Cauchy transform. In the
attempt to develop a boolean probability theory with amalgamation
over an algebra, when scalar-valued functionals are replaced by
conditional expectations, the above approaches are not suitable to
prove a similar multiplicative property, since objects such as the
Cauchy transform and resolvent function do not have natural
analogues. The present paper gives a combinatorial proof of
(\ref{multipli}), which is very easily adapted for the amalgamation
over an algebra case (see Section 4).

In \cite{dykema1}, \cite{dykema2}, K. Dykema develops an analogue
for the $S$-transform in the operator-valued free probability theory
and remarks that, due to non-commutativity, it satisfy a
``twisted'' multiplicative property. Namely, if $S=T^{-1}$, then\\
$T_{XY}=T_X(T_YIT_Y^{-1})T_Y$. Corollary \ref{finalcor} from Section
4 shows that the operator-valued $B$-transform has the same
multiplicative property as in the scalar case, without any
``twists''.

\section{Preliminaries}

Let $\cA$  and $\cC$ be a unital algebras and $\varphi:\cA\lra\cC$
be a linear functional with $\varphi(1)=1$. The  typical situation
is when $\cC$ is $\mathbb{C}$ or a subalgebra of $\cA$, but we
require it, nor that $\cC$ is commutative.

The subalgebras $\cA_1$ and $\cA_2$ of $\cA$ are said to be boolean
independent with respect to $\varphi$ (or just boolean independent)
if
\begin{equation}\label{1}
\varphi(X_1Y_1X_2\cdots)=\varphi(X_1)\varphi(Y_1)\cdots
\end{equation}
for all
 $X_1,X_2,\dots\in\cA_1$
  and
  $Y_1,Y_2,\dots\in\cA_2$.

If $\cA_1$ and $\cA_2$ are subalgebras of $\cA$, we will denote by
$\cA_1\bigvee\cA_2$ the algebra they generate in $\cA$. With this
notation, \ref{1} can be interpreted as follows: if $\cA_1$ and
$\cA_2$ are two boolean independent subalgebras of $\cA$, then
\begin{equation}\label{1b}
\varphi(a_1XYa_2)=\varphi(a_1X)\varphi(Ya_2)
\end{equation}
for all $X\in\cA_1, Y\in\cA_2$ and $a_1,a_2\in\cA_1\bigvee\cA_2$.

We define the boolean cumulant of order $n$ as the multilinear
function\\
 $b^n:\cA^n\lra\cC$ given by the recurrence
\begin{equation}
\label{2} \varphi(a_1\cdots a_n)= \sum_{k=1}^n b^k(a_1,\dots,a_k)
\varphi(a_{k+1}\cdots a_n).
\end{equation}
Note that (\ref{2}) is equivalent to
\begin{eqnarray*}
\varphi(a_1\cdots a_n) &=& \sum_{k=1}^n \varphi(a_1,\dots,a_k)
b^{n-k}(a_{k+1}\cdots a_n)\\
&=&\sum
b^{k_1}(a_1,\dots,a_{k_1})b^{k_2-k_1}(a_{k_1+1},\dots,a_{k_2})\cdots
b^{n-k_m}(a_{k_m+1},\dots,a_n)
\end{eqnarray*}
where the last summation is done over all $0\leq m\leq n-1$ and all
the $m$-uples $1\leq k_1<\dots<k_m$.
We will write $b^n_X$ for $b^n(X,\dots,X)$. For the moment,
respectively boolean cumulant generating power series of $X$ we will
use the notations $M_X(z)$, respectively $B_X(z)$, i.e.
\begin{eqnarray*}
M_X(z)&=&\sum_{n=1}^{\infty}\varphi(X^n)z^{n-1}\\
B_X(z)&=&\sum_{n=1}^{\infty}b_X^nz^{n-1}.
\end{eqnarray*}

The computations in the following two sections will often involve
functions of many arguments. For brevity sake, we will use the
short-hand notation described bellow:

If $I$ is a set of indices, $\{a_i\}_{i\in I}$ is a family of
elements from $\cA$ and $\ou=(i_1,\dots, i_n)$ is an ordered
sequence from $I$, we will write $\varphi(a_{\ou})$ for
$\varphi(a_{i_1}a_{i_2}\cdots a_{i_n})$, respectively $b^n(a_{\ou})$
for $b^n(a_{i_1},a_{i_2},\dots ,a_{i_n})$.  The length (cardinality)
of $\ou$ will be denoted by $|\ou|$. If $\ou_1,\ou_2$ are two
ordered sequences from $I$, their ordered concatenation will be
denoted by $\overrightarrow{u_1,u_2}$.

\section{Properties of the boolean cumulants}

The following property is analogue to the vanishing of free cumulant
with free independent entries. For the boolean case, only a weaker
property is true.

\begin{prop}\label{prop1}
Let $\cA_1,\cA_2$ be boolean independent subalgebras of $\cA$.
If $n,m\geq0$ and $a_1,\dots, a_{n+m+2}\in\cA_1\bigvee\cA_2$ such
that $a_{n+1}=X\in\cA_1$, $a_{n+2}=Y\in\cA_2$, then:
\[b^{n+m+2}(a_1,\dots,a_n, X,Y,a_{n+3},\dots, a_{n+m+2})=0\]
\end{prop}
\begin{proof}
Let $\ow=(1,2,\dots,n+m=2), \ou=(1,\dots,n), \ou_0=(1,\dots,n+1),
\ov=(n+3,\dots,n+m+2)$ and $\ov_0=(n=2,\dots, n+m=2$.

The proof will be done by induction on $n+m$. For $n=m=0$, (\ref{1})
implies
\[\varphi(XY)=\varphi(X)\varphi(Y)=b^1(X)\varphi(Y)\]
also, (\ref{2}) gives
\[\varphi(XY)=b^2(XY)+b^1(X)\varphi(Y),\]
therefore $b^2(XY)=0$.

Suppose the assertion is true for $n+m\leq N$ and let us prove it
for $n+m=N$. From (\ref{2}) we have that
\begin{eqnarray*}
\varphi(a_{\ow})&=&
\sum_{\ow=\overrightarrow{w_1,w_2}}
b^{|\ow_1|}(a_{\ow_1})\varphi(a_{\ow_2})\\
&=&
\sum_{
 \substack{\ow=\overrightarrow{w_1,w_2}\\ |\ow_1|\leq n+1}
}
b^{|\ow_1|}(a_{\ow_1})\varphi(a_{\ow_2})
 +
\sum_{
 \substack{\ow=\overrightarrow{w_1,w_2}\\ |\ow_1|> n+1}
  }
b^{|\ow_1|}(a_{\ow_1})\varphi(a_{\ow_2})
\end{eqnarray*}

On the other hand, from (\ref{1b}),
\begin{eqnarray*}
\varphi(a_{\ow})
&=&
\varphi(a_{\ou}XYa_{\ov})\\
&=&
\varphi(a_{\ou}X)\varphi(Ya_{\ov})
=\varphi(a_{\ou_0})\varphi(a_{\ov_0})\\
&=&
\sum_{\ou_0=\overrightarrow{u_1,u_2}}
b(a_{\ou_1})\varphi(a_{\ou_2})\varphi(a_{\ov_0})
\end{eqnarray*}
Since the last element in $a_{\ou_1}$ is $X\in\cA_1$ and the first
element in $a_{\ov_0}$ is $Y\in\cA_2$, property ($\ref{1b}$) gives
\[\varphi(a_{\ou_2})\varphi(a_{\ov_0})=\varphi(a_{\overrightarrow{u_2,v_0}}).\]
Therefore
\begin{eqnarray*}
\varphi(a_{\ow})
&=&
 \sum_{\ou_0=\overrightarrow{u_1,u_2}}
b(a_{\ou_1})\varphi(a_{\overrightarrow{u_2,v_0}})\\
&=&
\sum_
{
 \substack{\ow=\overrightarrow{w_1,w_2}\\ |\ow_1|\leq n+1}
}
b^{|\ow_1|}(a_{\ow_1})\varphi(a_{\ow_2})
\end{eqnarray*}
It follows that
\[
\sum_{
 \substack{\ow=\overrightarrow{w_1,w_2}\\ |\ow_1|> n+1}
  }
b^{|\ow_1|}(a_{\ow_1})\varphi(a_{\ow_2})
=0
\]
And the induction hypothesis implies the conclusion.
\end{proof}


\begin{cor}\label{additivity}
 If $X,Y$ are boolean independent, then
\[B_{X+Y}(z)=B_X(z)+B_Y(z).\]
\end{cor}

\begin{proof}
\ref{additivity} is implied by $b^n_{X+Y}=b^n_X+b^n_Y$, which is an
immediate consequence of \ref{prop1} and the multilinearity of the
mappings $b^n$.
\end{proof}

 The next proposition investigates properties of boolean cumulants with scalars
 among their entries.

\begin{prop}\label{prop2} If $n\geq1$ and $a_1,\dots,a_n\in\cA$, then
\begin{enumerate}
\item[(i)]$b^{n+1}(1,a_1,\dots,a_n)=0$
\item[(ii)]$b^{n+1}(a_1,\dots,a_n,1)=0$
\item[(iii)]$b^{n+1}(a_1,\dots,X_k,1,X_{k+1},\dots,a_n)=b^n(a_1,\dots,a_n).$
\end{enumerate}
\end{prop}
\begin{proof}

 \noindent(i): Let $\ou=(1,\dots,n),
\ov=(1,\dots,k),\ow=(k+1,\dots,n)$.
\[
\varphi(1,a_{\ou})=b^1(1)\varphi(a_{\ou}) +
\sum_{\ou=\overrightarrow{u_1,u_2}}b^{|\ou_1|+1}(1,a_{\ou_1})\varphi(a_{\ou_2})
\]
But $\varphi(1a_{\ou})=\varphi(a_{\ou})=b^1(1)\varphi(a_{\ou})$,
therefore
\[
\sum_{\ou=\overrightarrow{u_1,u_2}}b^{|\ou_1|+1}
(1,a_{\ou_1})\varphi(a_{\ou_2})=0.
\]
Since $\varphi(1X)=b^2(1,X)+b^1(1)\varphi(X)$, so $b^2(1,X)=0$,
property (ii) is proved by induction on $n$. The proof of (ii) is
similar.

\noindent(iii):From (i), one has that
\begin{eqnarray*}
\varphi(a_11a_2)
&=&
b^3(a-1,1,a_2)+b^2(a_1,1)\varphi(a_2)
+b^1(a_1)\varphi(1a_2)\\
&=& b^3(a_1,1,a_2)+b^1(a_1)\varphi(a_2).
\end{eqnarray*}
And, since $\varphi(a_1a_2)=b^2(a_1,a_2)+b^1(a_1)\varphi(a_2)$ we
have that $b^3((a_1,1,a_2)=b^2(a_1,a_2)$.

The rest of the proof is again a induction on $n$:
\begin{eqnarray*}
\varphi(a_{\ov}1a_{\ow})
&=&
 \sum_{\ov=\overrightarrow{v_1,v_2}}
 b^{|\ov_1|}(a_{\ov_1})\varphi(a_{\ov_2}1a_{\ow})
 +b^{|\ov|+1}\varphi(a_{\ow})+\\
 &&\sum_{\ow=\overrightarrow{w_1,w_2}}
 b^{k+1+|\ow_1|}(a_{\ou},1,a_{\ow_1}\varphi(a_{\ow_2})+
 b^{n+1}(a_{\ov},1,a_{\ou})
\end{eqnarray*}
Applying (ii) and the induction hypothesis, the previous equation
becomes
\begin{eqnarray*}
\varphi(a_{\ov}1a_{\ow})
 &=&
 \sum_{\ov=\overrightarrow{v_1,v_2}}
 b^{|\ov_1|}(a_{\ov_1})\varphi(a_{\ov_2}1a_{\ow})+
 \sum_{\ow=\overrightarrow{w_1,w_2}}
 b^{k+|\ow_1|}(a_{\ou},a_{\ow_1}\varphi(a_{\ow_2})\\
 &&
 +b^{n+1}(a_{\ov},1,a_{\ou})
\end{eqnarray*}
and (\ref{1}) gives the conclusion.
\end{proof}

\begin{cor}\label{cor2}
For all $X\in\cA$,$n\geq 2$, one has that

\[b^1_{1+X}=1+b^1_X\ \text{and}\ b^n_{X+1}=\sum_{k=0}^{n-2}{{n-2}\choose {k}}b_X^{k+2}.\]
\end{cor}
\begin{proof}
The multilinearity of $b^n$ implies that:
\[
b^n_{(1+X)}=\sum_{\substack{a_j\in\{1,X}\\j=1,\dots,n}b^n(a_1,a_2,\dots,a_n)
\]
 Terms having $1$ as first
or last entry also cancel from \ref{prop2}(i) and \ref{prop2}(ii).
The remaining terms can have the entry $1$ only on the remaining
$n-2$ positions left, and \ref{prop2}(iii) gives the stated result.
\end{proof}

Finally we will direct our attention toward boolean cumulants with
products among their entries. The multiplicative property of the
boolean cumulants will appear as a straightforward consequence of
the following proposition.

\begin{prop}\label{prop3}
Let $\cA_1$ and $\cA_2$ be two boolean independent subalgebras of
$\cA$, $X\in\cA_1, Y\in\cA_2$ and $a_1,\dots,a_{n}\in
\cA_1\bigvee\cA_2$. Then, for all $1\leq k\leq n$
\begin{enumerate}
\item[(i)] $b^{n+1}(a_{1},\dots,a_{k},XY,a_{k+1},\dots,a_n) =
b^{k+1}(a_{1},\dots,a_{k},X)b^{n-k+1}(Y,a_{k+1},\dots,a_n)$
\item[(ii)]$b^{n+2}(a_{1},\dots,a_{k},XY,XY,a_{k+1},\dots,a_n)=0.$
\item[(iii)]$b^{n+2}(a_{1},\dots,a_{k},Y,XY,a_{k+1},\dots,a_n)=0.$
\item[(iv)]$b^{n+2}(a_{1},\dots,a_{k},XY,X,a_{k+1},\dots,a_n)=0.$
\end{enumerate}
\end{prop}
\begin{proof}
Again we use the short-hand notations $\ou=(1,\dots,k),
\ov=(k+1,\dots,n)$.

\noindent(i): For $n=0$, $\varphi(XY)=b^1(XY)$; also
$\varphi(XY)=\varphi(X\cdot Y)=b^2(X,Y)+b^1(X)b^1(Y)$.
 But $b^2(X,Y)=0$ from \ref{prop1}, hence $b^1(XY)=b^1(X)b^1(Y)$.

Suppose now \ref{prop3}(i) true for $n\leq N$ and let us prove if
for $n=N+1$. One has that
\begin{eqnarray*}
\varphi(a_{\ou}XYa_{\ov})
 &=&
\varphi(a_{\ou}X)\varphi(Ya_{\ov})\\
&=&
 \sum_{\ou=\overrightarrow{u_1,u2}}
b^{|\ou_1|}(a_{\ou_1})\varphi(a_{\ou_2})\varphi(Ya_{\ov})
+b^{k+1}(a_{\ou},X)\varphi(Ya_{\ov})\\
 &=&
\sum_{\ou=\overrightarrow{u_1,u_2}}
b^{|\ou_1|}(a_{\ou_1})\varphi(a_{\ou_2})\varphi(Ya_{\ov})
+b^{k+1}(a_{\ou},X)b^1(Y)\varphi(a_{\ov})\\
&&
 +\sum_{\ov=\overrightarrow{v_1,v_2}}
  b^{k+1}(a_{\ou},X)b^{|\ov_1|+1}(Y,a_{\ov_1})\varphi(a_{\ov_2})\\
&&
  + b^{k+1}(a_{\ou},X)b^{n-k+1}(Y,a_{\ov}).
\end{eqnarray*}

Applying the induction hypothesis, the previous equality becomes
\begin{eqnarray*}
\varphi(a_{\ou}XYa_{\ov})
 &=&
 \sum_{\ou=\overrightarrow{u_1,u_2}}
b^{|\ou_1|}(a_{\ou_1})\varphi(a_{\ou_2})\varphi(Ya_{\ov})
+b^{k+1}(a_{\ou},XY)\varphi(a_{\ov})\\
&&
 +\sum_{\ov=\overrightarrow{v_1,v_2}}
  b^{k+|\ov_1|+1}(a_{\ou},XY,a_{\ov_1})\varphi(a_{\ov_2})\\
&&
  + b^{k+1}(a_{\ou},X)b^{n-k+1}(Y,a_{\ov}).
\end{eqnarray*}

Finally, comparing the above relation with the decomposition
(\ref{1b}) for\\ $\varphi(a_{\ou}\cdot XY \cdot a_{\ov})$, we obtain
\[
b^{n+1}(a_{\ou},XY,a_{\ov})=b^{k+1}(a_{\ou},X)b^{n-k+1}(Y,a_{\ov}).
\]
(ii):Part (i) gives
\begin{eqnarray*}
b^{|\ou|+|\ov|+2}(a_{\ou},XY,XY,a_{\ov})
 &=&
b^{|\ou|+1}(a_{\ou},X)\cdot b^{|\ov|+2}(Y,XY,a_{\ov})\\
&=&
b^{|\ou|+1}(a_{\ou},X)\cdot b^2(Y,X)\cdot b^{|\ov|+1}(Y,a_{\ov})\\
&=&0.
\end{eqnarray*}

(iii) and (iv) are analogous.
\end{proof}


\begin{cor}\label{c}For $\cA_1,\cA_2$ boolean independent algebras,
 $X\in\cA_1,Y\in\cA_2$ one has that:
\[b^n_{X+Y+XY}=b^n_X+b^n_Y+\sum_{k=1}^nb_X^kb^{n-k+1}_Y\]
\end{cor}
\begin{proof}
Let $\ou=(1,\dots,n)$. From the multilinearity of $b^n$,
\[b^n_{X+Y+XY}=\sum_{a_j\in\{X,Y,XY\}}b^n(a_{\ou})\]
From \ref{prop1}, all the boolean cumulants containing both $X$ and
$Y$, but no $XY$ vanish. From \ref{prop3}(ii)-(iv), among the
boolean cumulants containing $XY$, only ones of the form
$b^n(X,\dots,X,XY,Y,\dots, Y)$ are nonzero. Therefore:
\[b^n_{X+Y+XY}=
b^n_X+b^n_Y+\sum_{k=0}^{n-1}b^n(\underbrace{X,\dots,X}_{k\
\text{times}},XY,Y,\dots,Y)\] and the conclusion follows immediately
from \ref{prop3}(i).
\end{proof}

\begin{cor}
If $X,Y$ are boolean independent, then
\begin{equation}\label{multipli2}
B_{(1+X)(1+Y)}(z)=B_{(1+X)}(z)\cdot B_{(1+Y)}(z).
\end{equation}
\end{cor}
\begin{proof}
Let $Z=X+Y+XY$ and let $\alpha_n$, respectively $\beta_n$ be the
coefficients of $z^{n+1}$ from the left hand side, respectively the
right hand side on \ref{multipli2}. For $n=1$,
 \begin{eqnarray*}
 \alpha_1&=&b^1_{(1+X)(1+Y)}\\
 &=&
 1+b^1_X+b^1_Y+b^1_{XY}\\
 &=&
1+b^1_X+b^1_Y+b^1_{X}b^1_{Y}\\
&=& (1+b^1_X)(1+b^1_Y)=b^1_{(1+X)}b^1_{(1+Y)}\\
&=&\beta_1
 \end{eqnarray*}

 For $n=2$, we can write
 $\displaystyle{
 \alpha_n=\sum_{\substack{p,q\geq0\\p+q\geq1}}
 \alpha_{p,q}b^p_Xb^q_Y\ \text {and}\
 \beta_n=\sum_{\substack{p,q\geq0\\p+q\geq1}}
 \beta_{p,q}b^p_Xb^q_Y
 }$
 and (\ref{multipli2}) reduces to the equality between
 $\alpha_{p,q}$ and $\beta_{p,q}$.

Utilizing \ref{cor2},
\begin{eqnarray*}
\alpha_n
&=&
 b_{1+Z}^n\\
&=&
 \sum_{k=2}^n{{n-2}\choose {k-2}}b^k_Z\\
&=&
 \sum_{k=2}^n{{n-2}\choose {k-2}}
 \left(b^k_X+b^k_Y+\sum_{l=1}^kb_X^lb^{k+l+1}_Y\right)
\end{eqnarray*}
Therefore
\[
\alpha_{p,q}=
\begin{cases}
{{n-2}\choose {p+q-2}}& \text{if}\ p+q=n\ \text{and}\ p=0\ \text{or}\ q=0\\
{{n-2}\choose {p+q-3}}&\text{if}\ p,q\geq1 \ \text{and}\ p+q\leq n+1 \\
0& \text{otherwise}
\end{cases}
\]
On the other hand,
\begin{eqnarray*}
\beta_n
 &=&
 \sum_{k=1}^nb^k_{(1+X)}b^{n-k+1}_{(1+Y)}\\
 &=&
 b^1_{(1+X)}b^n_{(1+Y)}+b^n_{(1+X)}b^1_{(1+Y)}\\
 &&+\sum_{k=2}^{n-1}
 \left(
 \sum_{l=2}^k{{k-2}\choose{l-2}}b^l_X
 \sum_{s=2}^{n-k+1}{{n-k-1}\choose{s-2}}b_Y^{s}
 \right)
\end{eqnarray*}
Utilizing again \ref{cor2} and that $b^1_{1+X}=1+b^1_X$, we have
that
\begin{eqnarray*}
\beta_n
 &=&
 \sum_{k=2}^n{{n-2}\choose{k-2}}
 \left[
 b^k_X+b^k_Y+b^1_Xb^k_Y+b^k_Xb^1_Y
 \right]\\
 &&+\sum_{k=2}^{n-1}
 \left(
 \sum_{l=2}^k\sum_{s=2}^{n-k+1}{{k-2}\choose{l-2}}
 {{n-k-1}\choose{s-2}}b^l_Xb_Y^{s}
 \right)
\end{eqnarray*}
therefore
\[
\beta_{p,q}=
\begin{cases}
{{n-2}\choose {p+q-2}}
&
 \text{if}\ p+q=n\ \text{and}\ p=0\ \text{or}\ q=0\\
{{n-2}\choose {p+q-3}}
&
\text{if}\ p=1\ \text{or}\ q=1 \ \text{and}\ p+q\leq n+1 \\
\sum_{k=p}^{n-q+1}{{k-2}\choose{s-2}}{{n-k-1}\choose{q-2}}
 &
 \text{if}\ p,q\geq2,\ \text{and}\ p+q\leq n+1\\
 0& \text{otherwise}
\end{cases}
\]

The property (\ref{multipli2}) reduces though the to the equality
\[
\sum_{k=p}^{n-q+1}{{k-2}\choose{s-2}}{{n-k-1}\choose{q-2}}
={{n-2}\choose {p+q-3}}
\]
which is just an avatar of the well-known identity (see, for
example, \cite{GKP}, Chapter 5):
\begin{equation*}
\sum_{k=a}^{n-b}{{k}\choose{a}} {{n-k}\choose{b}}
={{n+1}\choose{a+b+1}}
\end{equation*}
\end{proof}


\section{Boolean independence with amalgamation over an algebra}
\subsection{Preliminaries}

 We need to consider an extended notion of
 non-unital complex algebra.
 $\gA$ will be
called a $\gB$-algebra
 if  $\gA$ is an algebra such that $\gB$ is a subalgebra of $\gA$
 or there is an algebra $\widetilde{\gA}$ containing $\gB$ as a
 subalgebra such that $\widetilde{\gA}=\gA\sqcup\gB$. (The symbol
 $\sqcup$ stands for disjoint union).

Let now $\gA$ be an algebra containing the unital subalgebra $\gB$
and $\Phi:\gA\lra\gB$ be a conditional expectation.

Suppose that $\gA_1$ and $\gA_2$ are two $\gB$-subalgebras of $\gA$.
 We say that $\gA_1$ and
$\gA_2$ are boolean independent with amalgamation over $\gB$ (or
just boolean independent over $\gB$) if
\[
\varphi(X_1Y_1X_2\cdots)=\varphi(X_1)\varphi(Y_1)\cdots
\]
for all
 $X_1,X_2,\dots\in\gA_1$
  and
  $Y_1,Y_2,\dots\in\gA_2$.

 Let $X\in\gA$ (if $\gA$ is a $\ast$-algebra, we also require
$X$ to be selfadjoint). If $X$ iS not commuting with $\gB$, the
algebra generated by $X$ and $\gB$ will not be spanned by
$\{\Phi(X^n)\}_n$. In this framework, the natural analogues of the
$n$-th moment
 of $X$ is the multilinear
function $\wm^n_X:{\gB}^{n-1}\lra\gB$ , given by
\[
\widetilde{m^n}_X(f_1,\dots,f_{n-1})=
\Phi(Xf_1X\cdots Xf_{n-1}X)
\]

The objects corresponding to power series are elements from $\mul$,
the set of multilinear function series over $\gB$ (see
\cite{dykema1}). A multilinear function series over $\gB$ is a
sequence $F=(F_0,F_1,\dots)$ such that $F_0\in\gB$ and $F_n$ is a
multilinear function from $\gB^n$ to $\gB$, for $n\geq1$. For
$F,G\in\mul$, the \emph{sum} $F+G$ and the \emph{formal product}
$FG$ are the elements from $\mul$ defined by:

\begin{eqnarray*}
(F+G)_n(f_1,\dots,f_n)&=&F_n(f_1,\dots,f_n)+G_n(f_1,\dots,f_n)\\
(FG)_n(f_1,\dots,f_n)&=&
\sum_{k=0}^nF_k(f_1,\dots,f_k)G_{n-k}(f_{k+1},\dots,f_n)\\
\end{eqnarray*}

\begin{defn}
The object corresponding to the $n$-th boolean cumulant $b^n$ (as
defined in Section 2) is the multilinear function
$\wb^n:\gA^n\times\gB^{n-1}\lra\gB$ given by the recurrence
\begin{eqnarray}
\Phi(X_1f_1X_2\cdots X_{n-1}f_{n-1}X_n)&=&\nonumber\\
&&\hspace{-2cm}\sum_{k=1}^n\wb^k_{X_1,\dots
X_k}(f_1,\dots,f_{k-1})f_k\Phi(X_{k+1}f_{k+1}\cdots
f_{n-1}X_n)\label{multirec1}
\end{eqnarray}
for all $X_1,\dots,X_n\in\gA$ and all $f_1,\dots f_{n-1}\in\gB$. The
$n$-uple $(X_1,\dots,X_n)$ will be called the \emph{lower argument}
of $\wb^n$, while $(f_1,\dots,f_{n-1}$ will be called the
\emph{upper argument} of $\wb^n$.

 If $X_1=\dots=X_n=X$,
 we will write $\wb^n_X$ for
$\wb^n_{X_1,\dots,X_n}$.
\end{defn}

We will denote by $\wM_X$, respectively $\wB$ the multilinear
function series for the multilinear moments and boolean cumulants of
$X\in\gA$. With this notation, the recurrence in the definition of
$\wb^n$ can be rewritten as
\begin{equation*}
\wM_X=\wB_X(1+I\wM_X)
\end{equation*}
 for $I$ the identity function
on $\gB$.
\begin{remark}\label{remark1}

Boolean independence over the algebra $\gB$ implies boolean
independence with respect to the functional $\Phi$, in the sense of
Section 2.

We denote by $b^n$ the $n$-th boolean cumulant with respect to
$\Phi$ in the sense of Section 2. Comparing the recurrences
(\ref{1b}) and (\ref{multirec1}), we have that

\begin{equation}\label{111}
b^n(X_1,\dots,X_n)=\wb^n_{X_1,\dots,X_n}(1,\dots,1).
\end{equation}
 Since, for all $f_0,h_1,\dots,f_{n-1}h_n\in\gB$,
\[
f_0\Phi(X_1h_1f_1X_2\cdots h_{n-1} f_{n-1}X_n)h_n=
\Phi((f_0X_1h_1)1(f_1X_2h_2)1\cdots 1(f_{n-1}X_nh_n))
\]
the recurrence (\ref{multirec1}) also gives
\begin{equation}\label{metamorf}
f_0[\wb^n_{X_1,\dots,X_n}(h_1f_1,\dots,h_{n-1}f_{n-1})]h_n
=b^n(f_0X_1h_1,\dots,f_{n-1}X_nh_n).
\end{equation}
\end{remark}

\subsection{Main results}

\begin{prop}\label{prop4.1}
 Let $\gA_1$ and $\gA_2$ independent over $\gB$ and
 $\gA_1\bigvee_\gB\gA_2$ be the subalgebra of $\gA$ generated by
 $\gA_1$,$\gA_2$ and $\gB$. Suppose that
$X\in\gA_1$, $Y\in\gA_2$ and $X_1\dots X_n\in\gA_1\bigvee_\gB\gA_2$.
 Then
 \begin{enumerate}
 \item[(i)]$\wb^{n+2}_{X_1,\dots,X_k, X,Y, X_{k+1},\dots X_n}=0$
\item[(ii)] for all $f\in\gB$, one has that $\displaystyle{
\wb^{n+1}_{f,X_1,\dots,X_n}=\wb^{n+1}_{X_1,\dots,X_n,f}=0,}$\
while\\
 ${}$\hspace{.1cm}$
\wb^{n+1}_{f,X_1,\dots,X_k,f,X_{k+1},\dots,X_n}(f_1,\dots,f_n)
=\wb^{n}_{X_1,\dots,X_n}(f_1,\dots,f_kff_{k+1},\dots,f_n) $
\item[(iii)]$\wb^{n+1}_{X_1,\dots,X_k,XY,X_{k+1},\dots,
X_n}(f_1,\dots,f_n)$\\
${}$\hspace{4cm} $=\wb^{k+1}_{X_1,\dots,X_k,X}(f_1,\dots,
f_k)\wb^{n-k+1}_{Y,X_{k+1},\dots,X_n,}(f_{k+1},\dots, f_n)$

\end{enumerate}
\end{prop}
\begin{proof}
The proposition is an immediate consequence of  Remark \ref{remark1}
and Propositions \ref{prop2} and \ref{prop3}.

For (i) we need to prove that
\[
\wb^{n+2}_{X_1,\dots,X_k, X,Y, X_{k+1},\dots
X_n}(f_1,\dots,f_{n+1})=0
\]
for all $f_1,\dots,f_{n+1}\in\gB$. But (\ref{metamorf}) implies
\[
\wb^{n+2}_{X_1,\dots,X_k, X,Y, X_{k+1},\dots X_n}(f_1,\dots,f_{n+1})
= b^{n+2}(X_1f_1,\dots,Xf_{k+1},Yf_{k+2},\dots,X_n)
\]
Since $Xf_{k+1}\in\gA_1$ and $Yf{k+2}\in\gA_2$, while
$X_1f_1,\dots,X_{n-1}f_{n+1},X_n\in\gA_1\bigvee_\gB\gA_2$,
Proposition \ref{prop1} implies the conclusion.

 The parts (ii) and (iii) are analogous consequences of
\ref{remark1} and \ref{prop2}(i), \ref{prop2}(ii), respectively
\ref{prop3}(i).
\end{proof}

For stating the next results, we need a brief discussion about
interval partitions.

\begin{defn}
An interval partition $\gamma$ on the set $\{1,2,\dots,n\}$ is a
collection of disjoint subsets, $D_1,\dots, D_q$, called
\emph{blocks}, such that
\begin{enumerate}
\item[(i)]$\bigcup_k=1^q D_k=\{1,\dots,n\}$
\item[(ii)]if $i_1<i_2$ and $k\in D_{i_1}$, $l\in D_{i_2}$, then
$k<l$.
\end{enumerate}
The set of all interval partitions on $\{1,\dots,n\}$ will be
denoted by $I(n)$. The number of blocks of the interval partition
$\gamma$ will be denoted by $|\gamma|$. For $\gamma_1\in I(n)$ and
$\gamma\in I(m)$, we will write $\gamma_1\oplus\gamma_2$ for the
interval partition from $I(m+n)$ obtained by juxtaposing $\gamma_1$
and $\gamma_2$.

If $\gamma\in\I(n-1)$,
$\gamma=(1,\dots,p_1)(p_1+1,\dots,p_2),\dots,(p_{k-1}+1,\dots,n-1)$,
and $(f_1,\dots,f_n)$ is a $n$-uple from $\gB$, and
 then $\pi_{\gamma}(f_1,\dots,f_n)$ will denote the
$q$-uple\\
$(f_1\cdots f_{p_1}, f_{p_1+1}\cdots f_{p_2},\dots,
f_{p_{q-1}+1}\cdots f_n)$ of elements from $\gB$.
\end{defn}
\begin{prop}
For all $X\in\gA$ and $n\geq2$, one has that
\begin{equation}\label{pp}
\wb^n_{(1+X)}(f_1,\dots,f_{n-1})= \sum_{\gamma\in\I(n-1)}
 \wb^{|\gamma|+1}_{X}
(\pi_\gamma(f_1,\dots,f_{n-1}))
\end{equation}
\end{prop}
\begin{proof}
From the multilinearity of $\wb^n$,
\begin{equation}\label{10}
\wb^n_{(1+X)}= \sum_{\substack{X_j\in\{1,X\}\\j=1,\dots,n}}
\wb^n_{X_1,\dots,X_n}=
 \sum_{\substack{X_j\in\{1,X\}\\j=2,\dots,n-1}}
\wb^n_{X,X_2,\dots,X_{n-1},X},
\end{equation}
since, from \ref{prop4.1}(ii), the terms with 1 on first and last
position vanish.

Let $\gamma\in\I(n-1)$,
$\gamma=(1,\dots,p_1)(p_1+1,\dots,p_2),\dots,(p_{k-1}+1,\dots,n-1)$.
Denote by $\wb^n_{\gamma, X}$ the boolean cumulant having the lower
argument $\underbrace{X,1,\dots,1}_{p_1},
\underbrace{X,1,\dots,1}_{p_2-p_1},\dots,$
$\dots,\underbrace{X,1,\dots,1}_{n-1-p_{k-1}},X$. Then (\ref{10})
becomes
\[
\wb^n_{(1+X)}(f_1,\dots,f_{n-1})
=\sum_{\gamma\in\I(n-1)}\wb^n_{\gamma,X}(f_1,\dots,f_{n-1})
\]
Finally, applying \ref{prop4.1}(iii), we get that
\begin{eqnarray*}
\wb^n_{\gamma,X}(f_1,\dots, f_{n-1}) &=&
\wb^{|\gamma|+1}_X(f_1f_2\cdots f_{p_1}, f_{p_1+1}\cdots
f_{p_2},\dots, f_{p_{k-1}+1}\cdots f_{n-1})\\
&=& \wb^{|\gamma|+1}_X(\pi_\gamma(f_1,\dots,f_n)),
\end{eqnarray*}
so q.e.d..
\end{proof}

\begin{cor}\label{finalcor}
Suppose $\gA_1,\gA_2$ are boolean independent over $\gB$ and
$X\in\gA_1$, $Y\in\gA_2$. Then
\[\wB_{X+Y}=\wB_X+\wB_Y\]
\[\wB_{(1+X)(1+Y)}=\wB_{(1+X)}\wB_{(1+Y)}.\]
\end{cor}
\begin{proof}
The first relation amounts to showing that
$\wb^n_{X+Y}=\wb^n_X+\wb^n_Y$, which is an obvious consequence of
Proposition \ref{prop4.1}(i) and of the multilinearity of $\wb^n$.

For the second relation, first note that the $n$-th component of
$B_{(1=x)(1+Y)}$ is
\begin{eqnarray*}
\wb^n_{(1+X)(1+Y)}(f_1,\dots,f_{n-1})
 &=&
  \wb^n_{1+(X+Y+XY)}(f_1,\dots,f_{n-1})\\
  &=&
  \sum_{\gamma\in\I(n-1)}
  \wb^{|\gamma|+1}_{X+Y+XY}(f_1,\dots,f_{n-1})
\end{eqnarray*}
The results \ref{prop4.1}(i) and (iii) imply the analogue of
Corollary \ref{c}:
\begin{eqnarray*}
\wb^m_{X+Y+XY}(f_1,\dots,f_{m-1})
 &=&
 \left(\wb^m_X(1+\wb^1_Y)+(1+\wb^1_X)\wb^m_Y\right)(f_1,\dots,f_{m-1})\\
 &&+\sum_{k=2}^{m-2}\wb^k_X(f_1,\dots,f_{k-1})\wb^{m-k}_Y(f_k,\dots,f_{m-k-1})
\end{eqnarray*}

Fix $\gamma\in\I(n-1)$. The above equation gives
\begin{eqnarray*}
\wb^{|\gamma|+1}_{X+Y+XY}(\pi_\gamma(f_1,\dots,f_{n-1}))
 &=&\left(\wb^{|\gamma|+1}_X(1+\wb^1_Y)+(1+\wb^1_X)\wb^{|\gamma|+1}_Y\right)
 (\pi_\gamma(f_1,\dots,f_{n-1}))\\
 &&\hspace{-2cm}+\sum_
  {
   \substack{
   \gamma_1\oplus\gamma_2=\gamma\\
\gamma_1\in\I(p)\\
\gamma_2\in\I(n-p-1)
 }
 }
 \wb^{|\gamma_1|+1}_X(\pi_{\gamma_1}(f_1,\dots,f_p))
 \wb^{|\gamma_2|+1}_y(\pi_{\gamma_2}(f_{p+1},\dots,f_{n-1}))
\end{eqnarray*}
therefore
\begin{eqnarray*}
\wb^n_{(1+X)(1+Y)}(f_1,\dots,f_{n-1})
 &=&
  \sum_{\gamma_\in\I(n-1)}\wb^{|\gamma|+1}_{\gamma,X+Y+XY}
  (\pi_\gamma(f_1,\dots,f_{n-1}))\\
 &&\hspace{-2.4cm}=
  \sum_{\gamma_\in\I(n-1)}
  \wb^{|\gamma|+1}_{X+Y+XY}(\pi_\gamma(f_1,\dots,f_{n-1}))\\
 &&\hspace{-2.4cm}=
  \sum_{\gamma\in I(n-1)}
  \left(
  \wb^{|\gamma|+1}_X(1+\wb^1_Y)+(1+\wb^1_X)\wb^{|\gamma|+1}_Y
   \right)
 (\pi_\gamma(f_1,\dots,f_{n-1}))\\
&&\hspace{-3.9cm}
  +\sum_{k=2}^{n-1}
  \left(
  \sum_{\gamma\in I(k-1)}\wb^{|\gamma|+1}_X(\pi_{\gamma}(f_1,\dots,f_{k-1}))
  \right)
   \left(
  \sum_{\sigma\in I(n-k)}\wb^{|\sigma|+1}_Y(\pi_{\sigma}(f_k,\dots,f_{n-k}))
  \right)
\end{eqnarray*}

Utilizing Prop(\ref{pp}), the above equation amounts to

\begin{eqnarray*}
\wb^n_{(1+X)(1+Y)}(f_1,\dots,f_{n-1})
 &=&
  \sum_{k=1}^n\wb^k_{(1+X)}(f_1,\dots,f_{k-1})\wb^{n-k+1}_{(1+Y)}(f_k,\dots,f_{n-1})
\end{eqnarray*}
hence the conclusion.

\end{proof}

 \textbf{Acknowledgements.}
My research was partially supported by the Grant 2-CEx06-11-34 of
the Romanian Government. I am thankful to \c{S}erban Belinschi, who
brought to my attention the problem addressed in \ref{finalcor}, to
Hari Bercovici for reference \cite{ber} and to Marek Bo\.{z}ejko for
 reference \cite{Schurmann}.

\end{document}